\documentclass{amsart}

\usepackage[T1]{fontenc}
\usepackage{lmodern}
\usepackage{amsmath,amssymb}
\usepackage{microtype}
\usepackage{hyperref}

\newtheorem{theorem}{Theorem}[section]
\newtheorem{lemma}[theorem]{Lemma}
\newtheorem{conjecture}[theorem]{Conjecture}
\numberwithin{equation}{section}
\newcommand{\R}{\mathbb{R}}
\newcommand{\Z}{\mathbb{Z}}

\title{Crescent configurations in three dimensions}
\author{Yuda Chen}
\email{math.yuda.chen@gmail.com}
\address{School of Mathematical Sciences and LPMC, Nankai University, Tianjin 300071, China}
\date{}

\begin{document}
	\begin{abstract}
		A crescent configuration of $n$ points in $\R^d$ is a set in general position whose distinct pairwise distances have multiplicities $1,\ldots,n-1$. We construct such a configuration in $\R^3$ for every $n\geq3$, giving a uniform bound of three for the minimum dimension. The proof uses equally spaced parameters on a short circular helix and Rolle's theorem. We also give explicit integer-coordinate configurations with $O(n\log n)$ bits per coordinate, and prove that every finite initial segment of the sequence $(\cos j,\sin j,2j)_{j\geq0}$ of length at least three is a crescent configuration.
	\end{abstract}
	\maketitle
	\section{Introduction}\label{sec:intro}
	For $d\geq2$, a finite set in $\R^d$ is in \emph{general position} if no $d+1$ of its points lie in an affine hyperplane and no $d+2$ lie on a sphere. Following Burt, Goldstein, Manski, Miller, Palsson, and Suh~\cite{BGMMPS}, an $n$-point set is a \emph{crescent configuration} if it is in general position and determines exactly $n-1$ distances, with multiplicities $1,\ldots,n-1$. Distances are Euclidean, and multiplicities count unordered pairs of distinct points; the distances need not be ordered by size. Without general position, equally spaced points on a line, or equally spaced angular parameters on a sufficiently short circular arc, give immediate examples for every $n$.
	The minimum dimension was studied by Burt et al.~\cite{BGMMPS}. Write $D(n)$ for the least dimension $d\geq2$ in which an $n$-point crescent configuration exists. They proved $D(n)\leq n-2$ for $n>3$, using a construction in which the ambient dimension increases with the number of points. Their argument builds configurations on spheres by adding points equidistant from the preceding points. In~\cite[Section~3]{BGMMPS}, they recorded a question of Albujer asking whether $D(n)$ remains bounded as $n\to\infty$. Our main result gives a uniform bound.
	\begin{theorem}
		\label{thm:main}
		For every integer $n\geq3$, the points
		\[
		P_j=\left(\cos\frac{j}{n},\ \sin\frac{j}{n},\ \frac{2j}{n}\right),
		\qquad j=0,\ldots,n-1,
		\]
		form a crescent configuration in $\R^3$. In particular,
		\[
		D(n)\leq3\qquad(n\geq3).
		\]
	\end{theorem}
	The points lie on a short circular helix. Their distance depends only on the difference of their indices, and different positive differences give different distances. The general-position conditions follow from Rolle's theorem: an arc with total angular variation less than $\pi$ meets each plane in at most three points and each sphere in at most four points. Thus the dimension bound has an elementary proof. It reduces the possible values of $D(n)$ to $2$ and $3$, but does not decide the planar existence problem.
	
	The planar problem goes back to Erd\H{o}s's Problem E2938 in the \emph{American Mathematical Monthly}~\cite{Erd82}; see also the later discussion~\cite{Erd89}. His 1986 account~\cite[p.~151]{Erd86} credits Pomerance with a five-point example, records a six-point example, and attributes independent seven-point constructions to Liu and Pal\'asti. These constructions appear in~\cite{Liu86,Pal87}. Pal\'asti subsequently constructed an eight-point example in the triangular lattice~\cite{Pal89}. Erd\H{o}s conjectured that planar examples do not exist for sufficiently large $n$~\cite{Erd86,Erd89}. Later computational work of Burt et al.~\cite[Remark~1]{BGMMPS} found no nine-point example among subsets of a specified $91$-point region of the triangular lattice; this was a search in a finite region, not a planar nonexistence theorem.
	
	Brass, Moser, and Pach~\cite[Section~5.10, Problem~2, p.~245]{BMP} include the planar existence question in their collection and record three-dimensional examples on up to six points due to Pal\'asti~\cite{PalastiSpace}. These finite examples precede the higher-dimensional construction of Burt et al.~\cite{BGMMPS}.
	
	Clemen, Dumitrescu, and Liu~\cite[Observation~5.1]{CDL} gave arbitrarily large planar sets with the same distance multiplicities that are not contained in a line or a circle. Their construction consists of $n-1$ points at equally spaced angular parameters on a circular arc of total angle less than $\pi/3$, together with the center of the circle. Since $n-1$ points are cocircular, these sets fail the general-position condition for $n\geq5$. Thus excluding a circle containing the entire set is weaker than the condition used here.
	In joint work with Durst, Hlavacek, Huynh, Miller, Palsson, and Wang~\cite{DHHMP}, the author developed an exact classification method based on complete distance-equality patterns up to vertex relabeling. Among the strict distance classes defined in that work, all three four-point classes and exactly $31$ of the $49$ five-point classes admit planar crescent realizations. The same work also reports that exactly six of the $6837$ six-point strict classes admit such realizations. The method separates finite combinatorial enumeration from geometric realizability, using Cayley--Menger and Gram-matrix criteria together with exact algebraic certificates of existence and nonexistence. Strictness is an additional combinatorial restriction, not a necessary condition for a planar crescent realization, so these results do not classify all distance classes admitting planar crescent realizations. The construction in Theorem~\ref{thm:main} instead realizes one distance-equality pattern in $\R^3$ uniformly for every $n$.
	
	Fish, King, Miller, Palsson, and Wahlenmayer~\cite{FKMPW} considered the problem in normed planes. They studied \emph{line-like} configurations, whose distance-equality pattern is that of equally spaced points on a line, and gave arbitrarily large examples for several classes of norms. They also introduced \emph{strong crescent configurations}, which satisfy an additional condition excluding four-point line-like subsets. That extra condition is not part of the higher-dimensional definition in~\cite{BGMMPS}. Our configurations are line-like, but satisfy both of its geometric general-position conditions.
	The same multiplicity requirement has also been studied without these Euclidean general-position restrictions. Demaine et al.~\cite{DemaineEtAl} classified the corresponding \emph{Erd\H{o}s-deep} sets in finite cyclic groups in their study of musical rhythms. Dukes and Gaede~\cite{DukesGaede} extended the problem to families of modular arithmetic progressions, including a classification for pairs of such progressions. These cyclic and normed-space variants concern different geometric realization problems from the one considered here.
	
	Helices themselves are familiar in the study of few-distance sets. Charalambides~\cite[Theorem~1.7]{Charalambides} showed that generalized helices support $n$-point sets with $O(n)$ distinct distances, obtained by taking equally spaced parameters. In our setting, the additional requirements are the exact multiplicities $1,\ldots,n-1$ and the two general-position conditions. A short helical arc meets all three requirements simultaneously.
	
	We give two further forms of the construction. First, the points can be chosen in the integer lattice, with an explicit bound on their coordinates.
	\begin{theorem}
		\label{thm:integer}
		For every integer $n\geq3$, there is an $n$-point crescent configuration in $\Z^3$ whose coordinates all have absolute value at most
		\[
		2(n-1)(4n^2+1)^{n-1}.
		\]
		In particular, each coordinate can be represented using $O(n\log n)$ bits.
	\end{theorem}
	Second, the angular step can be fixed independently of $n$.
	\begin{theorem}
		\label{thm:infinite}
		Let
		\[
		P_j=(\cos j,\sin j,2j),\qquad j=0,1,2,\ldots.
		\]
		No four points of this sequence are coplanar, and no five are cospherical. For every $n\geq3$, the set $\{P_0,\ldots,P_{n-1}\}$ is a crescent configuration in $\R^3$.
	\end{theorem}
	The integer-coordinate construction uses a rational rotation through a small angle. For the fixed infinite sequence, the short-arc argument no longer applies; instead, we use the transcendence of $e^i$ to establish the required determinant nonvanishing. The proofs of Theorems~\ref{thm:main}--\ref{thm:infinite} are given in Sections~\ref{sec:short}--\ref{sec:infinite}, respectively.
	\section{Short helical arcs}
	\label{sec:short}
	\begin{lemma}
		\label{lem:arc}
		Let $h\neq0$, and let $I\subset\R$ be an interval of length less than $\pi$. Among the points
		\[
		(\cos t,\sin t,ht),\qquad t\in I,
		\]
		no four are coplanar and no five are cospherical.
	\end{lemma}
	\begin{proof}
		Suppose a plane contains four distinct points of the arc. Substituting their coordinates into the equation of the plane gives four distinct zeros in $I$ of a function
		\[
		f(t)=a+bt+c\cos t+d\sin t,
		\qquad (b,c,d)\neq(0,0,0).
		\]
		
		By Rolle's theorem, $f''(t)=-c\cos t-d\sin t$ has at least two distinct zeros in $I$. A nonzero linear combination of $\cos t$ and $\sin t$ has consecutive zeros separated by $\pi$, so $c=d=0$. Then $b\neq0$, and the affine function $a+bt$ cannot have four zeros.
		
		If a sphere contains five distinct points of the arc, expanding its equation gives five distinct zeros in $I$ of a function
		\[
		g(t)=h^2t^2+a+bt+c\cos t+d\sin t.
		\]
		
		Three applications of Rolle's theorem show that $g'''(t)=c\sin t-d\cos t$ has at least two distinct zeros in $I$. Again $c=d=0$, whereas the remaining quadratic, with leading coefficient $h^2>0$, cannot have five distinct zeros.
	\end{proof}
	\begin{proof}[Proof of Theorem~\ref{thm:main}]
		The parameters $j/n$ lie in $[0,1]$, so Lemma~\ref{lem:arc}, with $h=2$, gives general position. If $k=|i-j|$, then
		\[
		\|P_i-P_j\|^2
		=2-2\cos\frac{k}{n}+4\left(\frac{k}{n}\right)^2.
		\]
		
		The expression $2-2\cos s+4s^2$ is strictly increasing for $0<s<1$, since its derivative is $2\sin s+8s>0$. Thus the gaps $k=1,\ldots,n-1$ give distinct distances. Exactly $n-k$ unordered pairs of indices have gap $k$, so the distance multiplicities are $n-1,n-2,\ldots,1$.
	\end{proof}
	\section{Integer coordinates}
	\label{sec:integer}
	\begin{proof}[Proof of Theorem~\ref{thm:integer}]
		Fix $n\geq3$, and choose the rational point on the unit circle
		\[
		z=\frac{4n^2-1+4ni}{4n^2+1}
		=e^{i\theta},
		\qquad \theta=2\arctan\frac{1}{2n}.
		\]
		Here $0<\theta<1/n$. The rational points
		\[
		P_j=(\Re(z^j),\Im(z^j),2j),\qquad 0\leq j<n,
		\]
		lie on the helix $(\cos t,\sin t,(2/\theta)t)$ at parameters $j\theta$. Since $(n-1)\theta<1$, Lemma~\ref{lem:arc} gives general position. Their squared distances are
		\[
		\|P_i-P_j\|^2=2-2\cos(k\theta)+4k^2,
		\qquad k=|i-j|.
		\]
		This is strictly increasing for $k=1,\ldots,n-1$: all the angles $k\theta$ lie in $(0,1)$, where cosine is strictly decreasing. As before, the distance associated with gap $k$ occurs $n-k$ times.
		Now multiply every coordinate by $(4n^2+1)^{n-1}$. Explicitly, take
		\begin{equation}
			\label{eq:integer-points}
			(4n^2+1)^{n-1}\bigl(\Re(z^j),\Im(z^j),2j\bigr),
			\qquad 0\leq j<n.
		\end{equation}
		The first two coordinates are integers because
		\[
		(4n^2+1)^{n-1}z^j
		=(4n^2+1)^{n-1-j}(4n^2-1+4ni)^j
		\]
		has integer real and imaginary parts. The third coordinate is also an integer. A common dilation preserves general position and distance multiplicities. Finally, $|z|=1$, so the first two coordinates in~\eqref{eq:integer-points} have absolute value at most $(4n^2+1)^{n-1}$, and the third is at most $2(n-1)(4n^2+1)^{n-1}$. The bit bound follows by taking logarithms.
	\end{proof}
	\section{A fixed infinite sequence}
	\label{sec:infinite}
	\begin{proof}[Proof of Theorem~\ref{thm:infinite}]
		For $P_j=(\cos j,\sin j,2j)$, the squared distance associated with an index gap $k$ is
		\[
		2-2\cos k+4k^2.
		\]
		Increasing $k$ to $k+1$ changes this expression by
		\[
		8k+4+2\bigl(\cos k-\cos(k+1)\bigr)\geq8k>0
		\qquad(k\geq1).
		\]
		Thus distinct positive gaps give distinct distances, and every initial segment has the required multiplicities.
		
		For general position, we use the Hermite--Lindemann theorem~\cite{Lindemann}. It implies that $e^i$ is transcendental, since $i$ is a nonzero algebraic number. Consequently, a nonzero Laurent polynomial with integer coefficients cannot vanish at $e^i$: multiplying by a sufficiently large power of the variable would otherwise give a nonzero integer polynomial having $e^i$ as a root.
		
		First fix integers $0\leq t_1<t_2<t_3<t_4$, and consider
		\begin{equation}
			\label{eq:plane-laurent}
			\det\begin{pmatrix}
				1&t_1&z^{t_1}&z^{-t_1}\\
				1&t_2&z^{t_2}&z^{-t_2}\\
				1&t_3&z^{t_3}&z^{-t_3}\\
				1&t_4&z^{t_4}&z^{-t_4}
			\end{pmatrix}.
		\end{equation}
		In the determinant expansion, the largest exponent of $z$ is $t_4-t_1$. It can arise only by selecting $z^{t_4}$ in the third column and $z^{-t_1}$ in the fourth. Its coefficient is therefore
		\[
		-\det\begin{pmatrix}1&t_2\\1&t_3\end{pmatrix}
		=-(t_3-t_2)\neq0.
		\]
		This proves that~\eqref{eq:plane-laurent} is not the zero Laurent polynomial, so it is nonzero at $z=e^i$. Replacing its last two columns by their half-sum and their difference divided by $2i$ is an invertible column transformation. It follows that the columns $1,t_r,\cos t_r,\sin t_r$, evaluated at these four indices, are linearly independent. Coplanarity of the corresponding points would give a nontrivial dependence among these columns, a contradiction.
		
		Next fix integers $0\leq t_1<\cdots<t_5$. The same argument applies to
		\begin{equation}
			\label{eq:sphere-laurent}
			\det\begin{pmatrix}
				1&t_1&t_1^2&z^{t_1}&z^{-t_1}\\
				1&t_2&t_2^2&z^{t_2}&z^{-t_2}\\
				1&t_3&t_3^2&z^{t_3}&z^{-t_3}\\
				1&t_4&t_4^2&z^{t_4}&z^{-t_4}\\
				1&t_5&t_5^2&z^{t_5}&z^{-t_5}
			\end{pmatrix}.
		\end{equation}
		The coefficient of its largest power, $z^{t_5-t_1}$, is
		\[
		\det\begin{pmatrix}
			1&t_2&t_2^2\\
			1&t_3&t_3^2\\
			1&t_4&t_4^2
		\end{pmatrix}
		=(t_3-t_2)(t_4-t_2)(t_4-t_3)\neq0.
		\]
		Thus~\eqref{eq:sphere-laurent} is nonzero at $e^i$, and the same column transformation shows that $1,t_r,t_r^2,\cos t_r,\sin t_r$ are linearly independent on these five indices.
		If the five points belonged to a sphere with center $(a,b,c)$ and radius $R$, expanding the sphere equation would give
		\[
		4t_r^2-4ct_r-2a\cos t_r-2b\sin t_r
		+1+a^2+b^2+c^2-R^2=0,
		\qquad r=1,\ldots,5.
		\]
		This is a nontrivial dependence among the five columns, since the coefficient of $t_r^2$ is $4$. The contradiction excludes cosphericity and completes the proof.
	\end{proof}
	\section{Final remark}
		In the preceding discussion, we have focused exclusively on $D(n)$. It should be noted that, as mentioned in Section~\ref{sec:intro}, the construction presented in this paper does not satisfy the criteria for a "strong crescent configuration" as defined in~\cite{FKMPW}. We conclude with a conjecture concerning strong crescent configurations. Let \(D_{\mathrm{strong}}(n)\) denote the minimum dimension in which an \(n\)-point strong crescent configuration exists.		
		\begin{conjecture}\label{conj:strong-crescent}
			We have
			\[
			D_{\mathrm{strong}}(n)=\Theta(n).
			\]
		\end{conjecture}
	\section*{Acknowledgments}
		The author thanks Prof.~Eyvindur Palsson for his helpful suggestions. OpenAI's GPT-5.5 Pro assisted in the search for the constructions presented in this paper. The author takes full responsibility for the mathematical arguments and results.
	\bibliographystyle{plain}
	\bibliography{bib}
\end{document}